\documentclass[english]{beamer}

\usepackage{amssymb}
\usepackage[cp1251]{inputenc}
\usepackage[english]{babel}
\usepackage{amsmath}
\usepackage{color}
\usepackage{xypic}
\usepackage{subfigure}
\usepackage{multicol}
\usepackage{epsfig}
\usepackage{graphicx}
\usepackage{url}
\usepackage{multimedia}
\usepackage{hyperref}

\def\mapr#1{\smash{\mathop{\buildrel{#1}\over\longrightarrow}}}

\def\C{{\bf C}}

\def\0{{\bf 0}}
\def\1{{\bf 1}}

\def\cA{{\cal A}}
\def\cB{{\cal B}}
\def\cC{{\cal C}}
\def\cD{{\cal D}}

\def\cH{{\cal H}}

\def\cL{{\cal L}}

\def\cT{{\cal T}}

\def\rg{{\frak g}}

\def\blf{\color{blue} $}
\def\elf{$ \color{black}}

\def\bdf{\color{blue}$$}
\def\edf{$$\color{black}}

\def\beq#1{\color{blue}\begin{equation}\label{#1}}
\def\eeq{\end{equation}\color{black}}

\def\mat{\color{blue} $}
\def\tam{$ \color{black}}

\def\mf{\color{blue} $}
\def\fm{$ \color{black}}

\def\Aut{{\hbox{\bf Aut}\;}}

\def\bydef{\stackrel{def}{=}}

\def\id{{\hbox{\bf Id}}}

\def\ker{{\hbox{\bf Ker}\;}}

\newtheorem{proposition}{Proposition}

\mode<presentation>

\AtBeginSection[]
{
  \begin{frame}<beamer>{Outline}
    \tableofcontents[currentsection]
  \end{frame}
}

\date{18.03.2011}

\title{{Characteristic classes of transitive Lie algebroids. Categorical point of view}}

\author{A.S.Mishchenko \\ (Harbin Institute of Technology, \\ Moscow State University)\\
Talk at the
International Conference\\
"Analysis, Topology and Applications",\\ Harbin, China}
\date{23.08.2011}

\begin{document}
\color{black}
\maketitle
\tableofcontents
\section{Introduction}
\begin{frame}
\frametitle{Introduction}
\framesubtitle{}
Transitive Lie algebroids have specific properties that allow to look at the
transitive Lie algebroid as an element of the object of a homotopy functor.
Roughly speaking each transitive Lie algebroids can be described as a vector
bundle over the tangent bundle of the manifold which is endowed with
additional structures.

Therefore transitive Lie algebroids admits a
construction of inverse image generated by a smooth mapping of smooth
manifolds.
\end{frame}

\begin{frame}
\frametitle{Introduction}
\framesubtitle{}
The construction can be managed as a homotopy functor
\mat\cT\cL\cA_{\rg}\tam from category of smooth manifolds
to the transitive Lie algebroids.

The functor \mat\cT\cL\cA_{\rg}\tam associates with
each smooth manifold \mat M\tam the set \mat\cT\cL\cA_{\rg}(M)\tam of all transitive algebroids with fixed structural finite dimensional Lie algebra \mat\rg\tam.

\end{frame}

\begin{frame}
\frametitle{Introduction}
\framesubtitle{}

The
intention of my talk is to use a homotopy classification of transitive Lie
algebroids due to K.Mackenzie

and on this basis to construct a classifying space.

\end{frame}

\begin{frame}
\frametitle{Introduction}
\framesubtitle{}

The
realization of the intention allows to describe characteristic classes of
transitive Lie algebroids form the point of view a natural transformation of
functors similar to the classical abstract characteristic classes for vector
bundles and to compare them with that derived from the Chern-Weil homomorphism by
J.Kubarski

\end{frame}
\section{Definitions and formulation of the problem}

\begin{frame}
\frametitle{Definitions}
\framesubtitle{}

Given smooth manifold \mat M\tam let

\bdf
E \mapr{a} TM \mapr{p_{T}} M
\edf

be a vector bundle over \mat TM\tam with fiber \mat\rg\tam. The fiber \mat\rg\tam has the structure of a finite dimensional Lie algebra
and the structural group of the bundle \mat E\tam is  \mat Aut(\rg)\tam, the group of all automorphisms of the Lie algebra \mat\rg\tam. Let
\mat p_{E}=p_{T}\cdot a\tam.
So we have a commutative diagram of two vector bundles
\bdf
\xymatrix{
E\ar[r]^{a}\ar[d]_{p_{E}}&TM\ar[d]^{p_{T}}\\
M\ar[r]&M
}
\edf

\end{frame}

\begin{frame}
\frametitle{Definitions}
\framesubtitle{}

The diagram is endowed with additional structure (commutator braces) and then is called
(Mackenzie, definition 3.3.1, Kubarski, definition 1.1.1)
transitive Lie algebroid
\bdf
\cA = \left\{
\makebox(70,25)[t]{\xymatrix{
E\ar[r]^{a}\ar[d]_{p_{E}}&TM\ar[d]^{p_{T}}\\
M\ar[r]&M
}}; \{\bullet,\bullet\}
\right\}.
\edf

\end{frame}

\begin{frame}
\frametitle{Definitions}
\framesubtitle{}

The braces \mat\{\bullet,\bullet\}\tam satisfy the natural properties, such that the space \mat\Gamma^{\infty}(E)\tam
with braces \mf\{\bullet,\bullet\}\fm forms an infinite dimensional Lie algebra with structure of the \mf C^{\infty}(M)\fm -- module,  that is

\end{frame}

\begin{frame}
\frametitle{Definitions}
\framesubtitle{}
that is
\begin{enumerate}
\item  Skew commutativity: for two smooth sections
    \mat\sigma_{1},\sigma_{2}\in \Gamma^{\infty}(E)\tam one has
    \beq{1}
    \{\sigma_{1},\sigma_{2}\}=-\{\sigma_{2},\sigma_{1}\}\in\Gamma^{\infty}(E),
    \eeq
\item Jacobi identity: for three smooth sections
    $\sigma_{1},\sigma_{2},\sigma_{3}\in \Gamma^{\infty}(E)$ one has
    \beq{3}
    \{\sigma_{1},\{\sigma_{2},\sigma_{3}\}\}+
    \{\sigma_{3},\{\sigma_{1},\sigma_{2}\}\}+
    \{\sigma_{2},\{\sigma_{3},\sigma_{1}\}\}=0,
    \eeq
\item Differentiation: for two smooth sections \mf\sigma_{1},\sigma_{2}\in \Gamma^{\infty}(E)\fm and smooth function \mf f\in C^{\infty}(M)\fm one has
    \beq{2}
    \{\sigma_{1},f\cdot\sigma_{2}\}=a(\sigma_{1})(f)\cdot\sigma_{2}+
    f\cdot\{\sigma_{1},\sigma_{2}\}
    \in\Gamma^{\infty}(E).
    \eeq

\end{enumerate}

\end{frame}

\begin{frame}
\frametitle{Pullback}
\framesubtitle{}

Let \mf f:M'\mapr{}M\fm be a smooth map. Then one can define an inverse image
(pullback) of the Lie algebroid (Mackenzie, page 156, Kubarski, definition 1.1.4),
\mf f^{!!}(\cA)\fm.

This means that given the finite dimensional Lie algebra \mf\rg\fm
there is the functor \mf\cT\cL\cA_{\rg}\fm such that
with any manifold \mf M\fm it assigns the family \mf\cT\cL\cA_{\rg}(M)\fm of all transitive  Lie algebroids
with fixed Lie algebra \mf\rg\fm.

\end{frame}

\begin{frame}
\frametitle{Pullback}
\framesubtitle{}

The following statement can be proved, see for example

\begin{theorem}
Each transitive Lie algebroid is locally trivial.
\end{theorem}

\end{frame}

\begin{frame}
\frametitle{Pullback}
\framesubtitle{}

This means that for a small neighborhood there is a trivia\-lization
of the vector bundles \mf E\fm, \mf TM\fm, \mf\ker a=L\approx \rg\times M\fm such that
\bdf
E\approx TM\oplus L,
\edf
and the Lie braces are defined by the formula:
\bdf
[(X,u),(Y,v)]= ([X,Y], [u,v]+X(v)-Y(u)).
\edf

\end{frame}

\begin{frame}
\frametitle{Homotopy of pullback}
\framesubtitle{}

Using the construction of pullback and the idea by Allen Hatcher

one can prove that the functor \mf\cT\cL\cA_{\rg}\fm is the homotopic functor.
More exactly for two homotopic smooth maps
\mf f_{0},f_{1}:M_{1}\mapr{}M_{2}\fm
and for the transitive  Lie algebroid
\bdf
\cA=\left(E\mapr{a}TM_{2}\mapr{}M_{2};\{\bullet,\bullet\}\right)
\edf
two inverse images \mf f_{0}^{!!}(\cA)\mf and \mf f_{1}^{!!}(\cA)\fm are isomorphic.

\end{frame}

\begin{frame}
\frametitle{Classifying space}
\framesubtitle{}
Let \mf\overline{\cT\cL\cA_{\rg}}\fm be the category of all transitive Lie algebroids and morphisms between them. The objects \mf\cA\in Obj(\overline{\cT\cL\cA_{\rg}})\fm are Lie algebroids
\bdf
\cA = \left\{
\makebox(70,25)[t]{\xymatrix{
E\ar[r]^{a}\ar[d]_{p_{E}}&TM\ar[d]^{p_{T}}\\
M\ar[r]&M
}}; \{\bullet,\bullet\}
\right\},
\edf
\mf M=M(\cA)\fm.

\end{frame}

\begin{frame}
\frametitle{Classifying space}
\framesubtitle{}

The morphisms \mf\varphi:\cA_{1}\mapr{}\cA_{2}\fm between the Lie algebroids are commutative diagrams
\bdf
{\xymatrix{
E_{1}\ar[r]^{\tilde f=E(\varphi)}\ar[d]_{p_{E_{1}}}&E_{2}\ar[d]^{p_{E_{2}}}\\
M_{1}\ar[r]_{f=M(\varphi)}&M_{2}
}}
\edf

\end{frame}

\begin{frame}
\frametitle{Classifying space}
\framesubtitle{}

Then one can define the direct limits
\bdf
\cB_{\rg}\bydef \lim\limits_{\mapr{}}\left(M(\cA);M(\varphi)\right); \cA\in Obj(\overline{\cT\cL\cA_{\rg}}),\varphi\in Morph(\overline{\cT\cL\cA_{\rg}})
\edf

\end{frame}

\begin{frame}
\frametitle{Classifying space}
\framesubtitle{}

Hence the final classifying space \mf\cB_{\rg}\fm has the property that
the family of all transitive  Lie algebroids with fixed Lie algebra \mf\rg\fm over the manifold \mf M\fm
has one-to-one correspondence with the family of homotopy classes of continuous maps \mf[M,\cB_{\rg}]\fm:
\bdf
\cT\cL\cA_{\rg}(M)\approx [M,\cB_{\rg}].
\edf

\end{frame}

\begin{frame}
\frametitle{Characteristic classes}
\framesubtitle{}

Using this observation one can describe the family of all characteristic classes of
a transitive  Lie algebroids in terms of cohomologies of the classifying space \mf\cB_{\rg}\fm. Really,
from the point of view of category theory a characteristic class \mf\alpha\fm is a natural transformation
from the functor \mf\cT\cL\cA_{\rg}\fm to the cohomology functor \mf\cH^{*}\fm.

\end{frame}

\begin{frame}
\frametitle{Characteristic classes}
\framesubtitle{}

This means that for the transitive  Lie algebroid
\mf\cA=\left(E\mapr{a}TM\mapr{}M;\{\bullet,\bullet\}\right)\fm the value of the characteristic class
\mf\alpha(\cA)\fm is a cohomology class
\bdf
\alpha(\cA)\in \cH^{*}(M),
\edf
such that for smooth map \mf f:M_{1}\mapr{}M\fm we have
\bdf
\alpha(f_{0}^{!!}(\cA))=f^{*}(\alpha(\cA))\in\cH^{*}(M_{1}).
\edf
\end{frame}

\begin{frame}
\frametitle{Characteristic classes}
\framesubtitle{}

Hence the family of all characteristic classes \mf\{\alpha\}\fm for  transitive
Lie algebroids with fixed Lie
algebra \mf\rg\fm has a one-to-one correspondence
with the cohomology group \mf\cH^{*}(\cB_{\rg})\fm.

On the base of these abstract considerations a natural problem can be
formulated.

\end{frame}

\begin{frame}
\frametitle{Problem}
\framesubtitle{}
\begin{problem}
Given finite dimensional Lie algebra \mf\rg\fm describe the
classifying space \mf\cB_{\rg}\fm
for  transitive Lie algebroids
in more or less understandable terms.
\end{problem}

Below we suggest a way of solution the problem and consider some trivial
examples.

\end{frame}
\section{Classifying cocycles due to K.Mackenzie}

\begin{frame}
\frametitle{Coupling}
\framesubtitle{}

Each transitive Lie algebroid

\bdf
\cA = \left\{
\makebox(70,25)[t]{\xymatrix{
E\ar[r]^{a}\ar[d]_{p^{E}}&TM\ar[d]^{p^{T}}\\
M\ar[r]&M
}}\quad , \{\bullet,\bullet\}
\right\}.
\edf
one can represent as an exact sequence of bundles

\bdf
0\mapr{}L\mapr{}E\mapr{a} TM\mapr{}0.
\edf
\end{frame}

\begin{frame}

\frametitle{Coupling}
\framesubtitle{}
The bundle \mat E \tam can be represent as a direct sum of bundles
\bdf
E=L\oplus TM,
\edf

\bdf
0\mapr{}L\mapr{}
\begin{array}{c}
L \\
\oplus \\
TM
\end{array}
\mapr{} TM\mapr{}0.
\edf
\noindent Then each section \mat \sigma\in\Gamma(E)\tam one can represent as the pair of sections
\bdf
\sigma = (u,X), \quad u\in\Gamma(L), \quad X\in\Gamma(TM).
\edf
\end{frame}

\begin{frame}
\frametitle{}
\framesubtitle{}

Then the commutator brace for the pair of the sections
\mat \sigma_{1}=(u_{1}, X_{1}), \quad \sigma_{2}=(u_{2}, X_{2})\tam
can be written by the formula

\bdf
\begin{array}{ll}
\{\sigma_{1},\sigma_{2}\}=&\{(u_{1}, X_{1}),(u_{2}, X_{2})\}=\\\\
&\hskip -1cm=\left(
[u_{1},u_{2}]+\nabla_{X_{1}}(u_{2})-\nabla_{X_{2}}(u_{1})+\Omega(X_{1},X_{2}),
[X_{1},X_{2}]
\right).
\end{array}
\edf

\end{frame}

\begin{frame}
\frametitle{}
\framesubtitle{}

Here
\bdf \nabla_{X}:\Gamma(L)\mapr{}\Gamma(L)\edf
is the covariant gradient of fiberwise differentiation of sections,
where
\bdf\Omega(X_{1},X_{2})\in \Gamma(L)\edf
is classical two-dimensional differential form with values in the fibers of the bundle  $L$.

The covariant derivative for fiberwise differentiation of the sections or so called linear connection in the bundle
$L$
\bdf \nabla_{X}:\Gamma(L)\mapr{}\Gamma(L)\edf
is an operator that satisfies the following natural conditions:

\end{frame}

\begin{frame}
\frametitle{}
\framesubtitle{}

\begin{enumerate}
\item Fiberwise differentiation with respect to multiplication in the Lie algebra structure of the fibre:
\bdf
\nabla_{X}([u_{1},u_{2}])=[\nabla_{X}(u_{1}),u_{2}]+
[u_{1},\nabla_{X}(u_{2})],\edf
\mat u_{1},u_{2}\in\Gamma(L)\tam.

\item Differentiation of sections in the space  \mat\Gamma(L)\tam as module over the function algebra  \mat \cC^{\infty}(M):\tam
\bdf
\nabla_{X}(f\cdot u)=X(f)\cdot u+f\cdot\nabla_{X}(u),
\edf
\mat u\in\Gamma(L)\tam, \mat f\in\cC^{\infty}(M)\tam.

\item Linear dependence on vector fields:
\bdf
\nabla_{f\cdot X_{1}+g\cdot X_{2}}=f\cdot \nabla_{X_{1}}+g\cdot \nabla_{X_{2}},
\edf
that is
\bdf
\nabla_{f\cdot X_{1}+g\cdot X_{2}}(u)=
f\cdot \nabla_{X_{1}}(u)+g\cdot \nabla_{X_{2}}(u),
\edf
\mat X_{1},X_{2}\in\Gamma(TM)\tam, \mat f,g\in\cC^{\infty}(M)\tam,
\mat u\in\Gamma(L)\tam.
\end{enumerate}

\end{frame}

\begin{frame}
\frametitle{}
\framesubtitle{}

More exactly for each vector field \mat X\in\Gamma(TM)\tam the covariant derivative associates a pair
\mat\nabla_{X}=(\cD,X)\tam,
\bdf
\cD:\Gamma(L)\mapr{}\Gamma(L), \quad
X:\cC^{\infty}(M)\mapr{}\cC^{\infty}(M),
\edf
that satisfies the conditions
\begin{enumerate}
\item \bdf
\cD([u_{1},u_{2}])=[\cD(u_{1}),u_{2}]+
[u_{1},cD(u_{2})],\quad u_{1},u_{2}\in\Gamma(L).\edf

\item
\bdf
\cD(f\cdot u)=X(f)\cdot u+f\cdot\cD(u),\quad u\in\Gamma(L),
\quad f\in\cC^{\infty}(M).
\edf
\end{enumerate}

\end{frame}

\begin{frame}
\frametitle{}
\framesubtitle{}

The association  \mat\nabla_{X}=(\cD,X)=(\cD(X),X)\tam satisfies the last condition
\begin{enumerate}
\addtocounter{enumi}{2}
\item \bdf
\nabla_{f\cdot X_{1}+g\cdot X_{2}}=f\cdot \nabla_{X_{1}}+g\cdot \nabla_{X_{2}},
\edf
\end{enumerate}

\end{frame}

\begin{frame}
\frametitle{}
\framesubtitle{}

The family of all covariant derivatives of fiberwise differentiation forms the infinite dimensional
Lie algebra with respect to operations of summation and composition. Really, let
\mat(\cD_{1},X_{1})\tam and \mat(\cD_{2},X_{2})\tam be two covariant derivatives.
The sum of derivatives
\mat(\cD_{3},X_{3})\tam is defined by the formula
\begin{enumerate}
\item \mat X_{3}=X_{1}+X_{2}\tam,
\item \mat\cD_{3}=\cD_{1}+\cD_{2}\tam
\end{enumerate}

The commutator brace
\mat(\cD_{3},X_{3})=\{(\cD_{1},X_{1}),(\cD_{2},X_{3})\} \tam
is defined by the formula:
\begin{enumerate}
\item \mat X_{3}=[X_{1},X_{2}]\tam,
\item \mat\cD_{3}=[\cD_{1},\cD_{2}]\tam
\end{enumerate}

\end{frame}

\begin{frame}
\frametitle{}
\framesubtitle{}

The family
of all covariant derivatives of fiberwise differentiation is
the space of sections of a transitive Lie algebroid, namely
\mat\cD_{der}(L)\mapr{}TM\tam, that is there is a bundle
\mat\cD_{der}(L)\mapr{}TM\tam, such that one has the exact sequence

\bdf
0\mapr{}Aut(L)\mapr{}\cD_{der}(L)\mapr{}TM\mapr{}0.
\edf
The bundle \mat\cD_{der}(L)\tam can be constructed as a union of fibres where each fiber
\mat\cD_{der}(L)_{x}\tam in the point \mat x\in M\tam consists of all covariant derivatives \mat(\cD,X)\tam in the point \mat x\in M\tam.
\end{frame}

\begin{frame}
\frametitle{}
\framesubtitle{}

That is the pair \mat(\cD,X)\tam
is the pair of operators
\bdf
\begin{array}{l}
\cD:\Gamma(L)\mapr{}L_{x}=\Gamma(x,L),
X:\C^{\infty}\mapr{}C,
\end{array}
\edf
which satisfy the conditions

\bdf
\begin{array}{l}
\cD([u_{1},u_{2}])=[\cD(u_{1}),u_{2}(x)]+
[u_{1}(x),\cD(u_{2})],\quad u_{1},u_{2}\in\Gamma(L),\\\\
\cD(f\cdot u)=X(f)\cdot u(x)+f(x)\cdot\cD(u),\\\\
u\in\Gamma(L)\quad f\in\cC^{\infty}(M).
\end{array}
\edf

\end{frame}

\begin{frame}
\frametitle{}
\framesubtitle{}

The fiber \mat\cD_{der}(L)_{x}\tam has finite dimension since belongs to the exact sequence
\bdf
0\mapr{}Aut(L_{x})\mapr{}\cD_{der}(L)_{x}\mapr{}T_{x}M\mapr{}0,
\edf
and the kernel  \mat Aut(L_{x}) \tam consist of the operators  of the form \mat(\cD,X)\tam,
that satisfy the conditions
\bdf
\begin{array}{l}
\cD([u_{1},u_{2}])=[\cD(u_{1}),u_{2}(x)]+
[u_{1}(x),\cD(u_{2})],\quad u_{1},u_{2}\in\Gamma(L),\\\\
\cD(f\cdot u)=f(x)\cdot\cD(u),\\\\
u\in\Gamma(L)\quad f\in\cC^{\infty}(M),
\end{array}
\edf
hence consists of automorphisms of the finitely dimensional Lie algebra \mat L_{x}\tam.

\end{frame}

\begin{frame}
\frametitle{}
\framesubtitle{}

The exact sequence of the bundles

\bdf
0\mapr{}Aut(L)\mapr{}\cD_{der}(L)\mapr{}TM\mapr{}0.
\edf

can be included in the exact diagram

\end{frame}

\begin{frame}
\frametitle{}
\framesubtitle{}

\beq{3}
\xymatrix{
&ZL\ar[r]^{=}\ar[d]^{i}&ZL\ar[d]^{i}\\
&L\ar[r]^{=}\ar[d]^{ad} &L\ar[d]^{ad}\\
0\ar[r]&Aut(L)\ar[r]^{j}\ar[d]^{\natural^{0}}&\cD_{Der}(L)\ar[r]^{a}\ar[d]^{\natural}&TM\ar[r]\ar[d]^{=}&0\\
0\ar[r]&Out(L)\ar[r]^{\bar j}\ar[d]&Out\cD_{Der}(L)\ar[r]^{\bar a}\ar[d]&TM\ar[r]\ar[d]&0\\
&0&0&0
}
\eeq

\end{frame}

\begin{frame}
\frametitle{}
\framesubtitle{}

The splitting of the algebroid \mat E\tam  in a direct sum
\bdf
E=L\oplus TM,
\edf
means that there is a splitting map in the exact sequence
\bdf
0\mapr{}L\mapr{}E\mapr{a} TM\mapr{}0,
\edf
\mat T:TM\mapr{} E\tam, \mat a\circ T=\id\tam,
that is called transversal.

\end{frame}

\begin{frame}
\frametitle{}
\framesubtitle{}

There is a morphism of the algebroids
\bdf
\xymatrix{
E\ar[rr]^{ad}\ar[dr]\ar[rdd]&&\cD_{der}(L)\ar[dl]\ar[ldd]\\
&TM\ar[d]&\\
&M&
}
\edf
so called adjoint representation of the Lie algebroid \mat E\tam.

\end{frame}

\begin{frame}
\frametitle{}
\framesubtitle{}

This representation also can be included in an exact diagram
\bdf
\xymatrix{
0\ar[r]&L\ar[r]\ar[d]^{ad}&E\ar[r]\ar[d]^{ad}&TM\ar[r]\ar[d]^{=}&0\\
0\ar[r]&Aut(L)\ar[r]&\cD_{der}(L)\ar[r]&TM\ar[r]&0
}
\edf

\end{frame}

\begin{frame}
\frametitle{}
\framesubtitle{}

and can be extended to the following diagram
\bdf
\xymatrix{
&0\ar[d]\\
&ZL\ar[d]\\
0\ar[r]&L\ar[r]\ar[d]^{ad}&E\ar[r]\ar[d]^{ad}&TM\ar[r]\ar[d]^{=}&0\\
0\ar[r]&Aut(L)\ar[r]\ar[d]^{\natural^{0}}&\cD_{der}(L)\ar[r]&TM\ar[r]&0\\
&Out(L)\ar[d]\\
&0
}
\edf

\end{frame}

\begin{frame}
\frametitle{}
\framesubtitle{}

and can be closed to the following diagram

\beq{4}
\xymatrix{
&ZL\ar[d]\ar[r]^{=}&ZL\ar[d]\\
0\ar[r]&L\ar[r]\ar[d]^{ad}&E\ar[r]\ar[d]^{ad}&TM\ar[r]\ar[d]^{=}&0\\
0\ar[r]&Aut(L)\ar[r]\ar[d]^{\natural^{0}}&\cD_{der}(L)\ar[r]\ar[d]^{\natural^{!}}&TM\ar[r]&0\\
&Out(L)\ar[d]\ar[r]^{=}&Out(L)\ar[d]\\
&0&0
}
\eeq

\end{frame}

\begin{frame}
\frametitle{}
\framesubtitle{}

Hence

\beq{5}
\xymatrix{
&ZL\ar[d]\ar[r]^{=}&ZL\ar[d]\\
0\ar[r]&L\ar[r]\ar[d]^{ad}\ar[dr]^{ad}&E\ar[ddr]^{\natural\circ ad}\ar[rr]^{a}\ar[d]^{ad}&&TM\ar@{.>}[ddl]_{T}\ar[r]\ar[d]^{=}&0\\
0\ar[r]&Aut(L)\ar[r]\ar[dd]^{\natural^{0}}&\cD_{der}(L)\ar[dr]^{\natural}\ar[rr]^{a}\ar[dd]^{\natural^{!}}&&TM\ar[r]&0\\
&&&Out\cD_{der}(L)\ar[ur]^{a}\ar[dr]\ar[dl]\\
&Out(L)\ar@{.>}[urr]\ar[r]^{=}&Out(L)&&0\\
}
\eeq

\end{frame}

\begin{frame}
\frametitle{}
\framesubtitle{}

All homomorphisms of the bundles are morphisms of the Lie algebroids. In particular the exact sequence

\beq{6}
\xymatrix{
0\ar[r]&Out(L)\ar[r]^{\bar j}&Out\cD_{Der}(L)\ar[r]^{\bar a}&TM\ar[r]&0\\
}
\eeq
has a transversal \mat T\tam which also is a morphism of the Lie algebroids. This means that
\mat T\tam is a coupling.
\bdf
\xymatrix{
0\ar[r]&Out(L)\ar[r]^{\bar j}&Out\cD_{Der}(L)\ar[r]_{\bar a}&TM\ar[r]\ar[l]_{T}&0.
}
\edf
Therefore the bundle \mat Out(L)\tam is flat.

\end{frame}

\begin{frame}
\frametitle{}
\framesubtitle{}

In (Mackenzie, Lemma 7.2.4) a cocycle \mf f_{T}\in \Omega^{3}(M,\rho^{T}, ZL)\fm was defined that plays the role of the obstruction \mf f_{T}\in Obs(T)\in \cH^{3}(M,\rho^{T},ZL)\fm for existence of the structure of the Lie algebroids \mf\cA\fm.
The key property of the cocycle \mf f_{T}\fm is that it is functorial. Precisely, if \mf\varphi:M_{1}\mapr{}M_{2}\fm, \mf L_{1}=\varphi^{*}(L_{2})\fm, \mf T_{1}\fm is the inverse image of the coupling \mf T_{2}\fm then
\bdf
\varphi^{*}(Obs(T_{2}))= Obs(T_{1}).
\edf

\end{frame}

\section{Classifying space}

\begin{frame}
\frametitle{}
\framesubtitle{}

The bundle $L$ has the structural group $Out(L)^{\delta}$ with more fine topology such that the exact sequence

\bdf
\xymatrix{
0\ar[r]&Z\rg\ar[r]&\rg\ar[r]^{ad}&\Aut(\rg)\ar[r]^{\natural^{0}}&Out(\rg)_{discrete}\ar[r]&0
}
\edf
has the discrete topology in the group \mf Out(\rg)\fm.

\end{frame}

\begin{frame}
\frametitle{}
\framesubtitle{}

Hence the bundle \mf L\fm can be described as inverse image from classifying space of the structural group \mf Out(\rg)^{\delta}\fm:
\bdf
\varphi: M\mapr{} BOut(\rg)^{\delta};
\edf
\blf
L=\varphi^{*}(L^{\infty})
\elf

\end{frame}

\begin{frame}
\frametitle{}
\framesubtitle{}

\begin{theorem}
There is a cohomology class
\mf Obs\in\cH^{3}(BOut(\rg)^{\delta}, \rho, Z\rg)\fm
such that
\bdf
Obs(L)=\varphi^{*}(Obs)\in \Omega^{3}(M,\rho^{T}, ZL).
\edf
\end{theorem}

\end{frame}

\begin{frame}
\frametitle{}
\framesubtitle{}

Due to

the obstruction class \mf Obs\in\cH^{3}(BOut(\rg)^{\delta}, \rho, Z\rg)\fm can be identify with homotopy class of continuous maps
\bdf
\xymatrix{
BOut(\rg)^{\delta}\ar[r]^{\varphi_{Obs}}&
K_{Out(\rg)}(Z\rg, 3)\\
}
\edf
where \mf K_{G}(Z, 3)\fm denotes the equivariant version of the
Eilinberg-MacLane complex for \mf G\fm--module \mf Z\fm.

\end{frame}

\begin{frame}
\frametitle{}
\framesubtitle{}
The Eilinberg-MacLane complex \mf K_{G}(Z, 3)\fm has the subcomplex \mf K(G,1)\subset K_{G}(Z, 3)\fm which pays the role of the point for classical case. Then we have the
diagram of fibrations
\bdf
\xymatrix{
BOut(\rg)^{\delta}\ar[r]^{\varphi_{Obs}}&
K_{Out(\rg)}(Z\rg, 3)\\
F(\rg)\ar[u]\ar[r]&K(Out(\rg),1)\ar[u]
}
\edf

\end{frame}

\begin{frame}
\frametitle{}
\framesubtitle{}
For a transitive Lie algebroid
 \bdf
\cA = \left\{
\makebox(70,25)[t]{\xymatrix{
E\ar[r]^{a}\ar[d]_{p^{E}}&TM\ar[d]^{p^{T}}\\
M\ar[r]&M
}}\quad , \{\bullet,\bullet\}
\right\}
\edf

 with fixed structural Lie algebra \mf\rg\fm one has a continuous map
\bdf
f: M\mapr{}F(\rg)
\edf
\end{frame}

\begin{frame}
\frametitle{}
\framesubtitle{}
that is generated by a continuous map
\bdf
\cB_{\rg}\mapr{}F(\rg)
\edf
That is one has the commutative diagram
\bdf
\xymatrix{
&[M,BOut(\rg)^{\delta}]\\
\cT\cL\cA_{\rg}(M)\approx [M,\cB_{\rg}]\ar[r]\ar[ur]^{Coupling(T)}&[M, F(\rg)]\ar[u]
}
\edf
\end{frame}

\begin{frame}
\frametitle{Classifying space}
\framesubtitle{}
Due to Mackenzie theorem [7.3.18]
\begin{theorem}
Let \mf\cA\fm be a Lie algebroid on \mf M\fm, let \mf L\fm be a Lie algebra bundle on \mf M\fm and
let \mf T\fm be a coupling of \mf\cA\fm with \mf L\fm such that
\bdf \mbox{Obs}(T)=0\in \cH^{3}(\cA, \rho^{T}, ZL).\edf
Then the additive group \mf\cH^{2}(\cA, \rho^{T}, ZL)\fm acts freely and transitively on the fiber of the map \mf Coupling(T)\fm.
\end{theorem}

\end{frame}

\begin{frame}
\frametitle{Classifying space. Characteristic classes}
\framesubtitle{}
Therefore
\bdf
\cB_{\rg}\approx F(\rg)\times K_{Out(\rg)}(Z\rg, 2)
\edf
\begin{proposition}
The family of all characteristic classes for transitive Lie
algebroids \mf \cT\cL\cA_{\rg}(M)\fm can be identified with
\mf \cH^{*}
(F(\rg)\times K_{Out(\rg)}(Z\rg, 2))\fm
\end{proposition}

\end{frame}
\section{Characteristic classes}

\begin{frame}
\frametitle{Chern-Weil homomorphism}
\framesubtitle{}
The Chern-Weil homomorphism is a homomorphism from the algebra of invariants
\mf\bigoplus\left(\Gamma\bigvee L^{*}\right)_{I}\fm
to the cohomology of the Lie algebroid \mf\cA\fm, \mf\cH^{*}(M,\cA)\fm:

\bdf
h_{\cA}:\bigoplus\left(\Gamma\bigvee L^{*}\right)_{I}
\mapr{}\cH^{*}(M,\cA)
\edf

that gives a nonclassical example of the characteristic classes since the Chern-Weil homomorphism commutes with pullback of the Lie algebroids (Kubarski, theorem 4.2.2, theorem 4.3.1)
\end{frame}

\begin{frame}
\frametitle{Chern-Weil homomorphism}
\framesubtitle{}
The difference consists on that the cohomology for Kubarski case depends on the choice of the Lie algebroid on \mf M\fm.
So we should introduce here the cohomology of the classifying space \mf\cB_{\rg}\fm with coefficients in nonexisting Lie algebroid on the \mf\cB_{\rg}\fm.
\end{frame}

\begin{frame}
\frametitle{Chern-Weil homomorphism}
\framesubtitle{}
\begin{block}{}
Partly supported by RFBR 11-01-00057-a, 10-01-92601-KO\_a, 11-01-90413-Ukr\_f\_a and State program RNP 2.1.1.5055
\end{block}
\end{frame}
\begin{frame}
\frametitle{References}
\framesubtitle{}

\end{frame}

\begin{frame}
\frametitle{References}
\framesubtitle{}

\end{frame}

\end{document}